\documentclass[copyright,creativecommons]{eptcs}
\usepackage[utf8]{inputenc}
\usepackage[T1]{fontenc}
\usepackage{hyperref}
\usepackage{amssymb}
\usepackage{tikz-cd}
\usepackage{mathtools}
\usepackage{prftree}

\providecommand{\event}{FICS 2024}

\usepackage{amsmath}
\usepackage{amssymb}
\usepackage{amsthm}

\theoremstyle{definition}
\newtheorem{thm}{Theorem}
\newtheorem{lem}[thm]{Lemma}
\newtheorem{df}[thm]{Definition}
\newtheorem{cor}[thm]{Corollary}

\newcommand{\construct}{\textsf{construct}}
\newcommand{\destruct}{\textsf{destruct}}
\newcommand{\subtree}[2]{\mathsf{ST}_{#2}(#1)}
\newcommand{\restricts}{{\upharpoonright}}

\newcommand{\sequent}{\text{Seq}}
\newcommand{\rul}{\text{Rul}}

\newcommand{\tree}{\mathbb{T}}
\newcommand{\ftree}{\mathbb{T}^{<\omega}}
\newcommand{\arity}[2]{#1\textsf{-arity}(#2)}
\newcommand{\pf}[1]{\mathbb{P}(#1)}

\ifcsname endofdump\endcsname\endofdump\fi

\author{Borja Sierra Miranda\footnote{This work is based on a master project done at the Master of Logic in the University of Amsterdam. I thank Bahareh Afshari and Lide Grotenhuis for supervising the project.}\\Logic and Theory Group, University of Bern\\borja.sierra@unibe.ch}
\date{}
\title{Cyclic Proofs for iGL via Corecursion}
\def\titlerunning{Cyclic Proofs for iGL via Corecursion}
\def\authorrunning{B. Sierra Miranda}

\begin{document}

\maketitle

\section{Introduction}
\label{sec:orgef226cc}

Cyclic proof theory studies proofs where cycles are allowed (see \cite{brotherston2006sequent}). This is useful for developing proof theory for logics with fixpoint operators: cycles can be used to represent the unfolding of a fixpoint.
However, this cyclic character is not unique to such explicit fixpoints. For example, modal logics whose frames have a Noetherian (conversely wellfounded) condition, such as \(\textsf{GL}\) (G\"odel-L\"ob logic, see \cite{ShamkanovGl},\cite{iemhoff2021reasoning}), \(\textsf{S4Grz}\) (Grzegorczyk logic, see \cite{ShamkanovGrz}) and \(\textsf{K4Grz}\) (see \cite{guillermo}) also have cyclic proof systems. 

Particularly, in \cite{ShamkanovGl}, Shamkanov introduces a non-wellfounded and a cyclic sequent system (\(\textsf{GL}\)). He proves the equivalence of these two systems with an acyclic finite system via proof translations. In order to go from the finite system to the non-wellfounded system he defines the translation by corecursion.

In \cite{iemhoff2021reasoning}, Iemhoff generalized the work of Shamkanov studying when, for a given modal logic proof system, there exists another modal logic proof system such that proofs in the first are equivalent to cyclic proofs in the second. There, she shows that \(\textup{i}\mathsf{GL}\), an\footnote{The use of ``an'' instead of ``the'' is deliberate. Check the footnote of page \pageref{X} for an explanation.} intuitionistic version of \(\textsf{GL}\), also has a natural cyclic proof system.

We provide an alternative proof of the equivalence of a standard calculus for \(\textup{i}\mathsf{GL}\) and a cyclic one. The difference from the above-mentioned previous work is twofold:
\begin{enumerate}
    \item Part of the motivation of \cite{iemhoff2021reasoning} is not to use a non-wellfounded system as in \cite{ShamkanovGl}. We want to use a non-wellfounded system and in particular define the proof translation by corecursion.
    \item In \cite{ShamkanovGl} and \cite{iemhoff2021reasoning} cyclic systems are defined such that every cycle is allowed. We are going to see that for the sequent calculus of \(\textup{i}\mathsf{GL}\) that we work with this is not possible. The difference is that for propositional logic we will use the rules of the intuitionistic calculus \(\textsf{G3}\text{i}\) (see \cite{Troelstra_Schwichtenberg_2000}).  As we have already mentioned, Iemhoff's method can be applied to obtain a cyclic calculus for \(\textup{i}\mathsf{GL}\), but using Dyckhoff calculus (see \cite{Dyckhoff}) instead of \(\textsf{G3}\text{i}\). See Section \ref{sec:progress} for details.
\end{enumerate}

\section{Algebras and Coalgebras: (co)recursion\protect\footnote{We are just going to introduce the necessary principles for our work, particularly infinite trees and corecursion to them. These ideas can already found in \cite{COURCELLE198395}, so the methods we use can be considered to be standard for the treatment of infinite trees. For a general introduction on coalgebras the reader can consult \cite{VENEMA2007331}.}}

We need to work with non-wellfounded trees. We are going to use the representation of trees with a non-empty set of finite sequences of natural numbers. Given a set \(A\) we will write \(A^*\) to denote the set of finite sequences over \(A\).  Elements of \(\omega^*\) will be denoted by \(w,v,u\).

\begin{df}
An \emph{\(L\)-labelled tree} is a pair \(T = (N,\ell)\), where:
\begin{enumerate}
\item \(N \subseteq \omega^*\), non-empty and closed under initial segments.
\item For any \(w \in N\) there exists a natural number \(m\) such that for any \(i\), \(wi \in N\) iff \(i < m\). Given \(w\), this number can be shown to be unique and we call it the \emph{arity of \(w\) in \(T\)}, \(\arity{T}{w}\).
\item \(\ell : N \longrightarrow L\), called the labelling function of \(T\).
\end{enumerate}
\(T\) is said to be \emph{finite} iff \(N\) is finite. A branch of \(T\) is just an infinite path of \(T\). We denote the collection of \(L\)-labelled trees as \(\tree_{L}\) and the collection of finite \(L\)-labelled trees as \(\ftree_{L}\). In case it can be filled by context, we will omit $L$.
\end{df}
Given a finite or infinite sequence \(w\) and a natural number \(i\), we will write \(w \restricts i\) to mean the restriction of \(w\) to \(\{0,\ldots,i-1\}\). 
If \(X\) is a set of finite sequences, we say that \(w\) is maximal in \(X\) iff there is no sequence in \(X\) strictly extending \(w\).
Let us define some usual notions of trees in this formalism.
\begin{df}
    A \emph{branch} of \(T\) is just an infinite sequence \(b \in \omega^\omega\) such that for any \(i \in \omega\), \(b{\upharpoonright} i\) is a node of \(T\).

    A \emph{leaf} of \(T\) is a maximal sequence in the nodes of \(T\). An \emph{internal node} is any node that is not a leaf.

    Let \(T=(N,\ell)\) be a tree and \(w\) one of its nodes. We define the \(T\)\emph{-subtree generated at }\(w\), as \(\subtree{T}{w} = (N',\ell')\) where:
    \begin{align*}
        &N' = \{v \in \omega^* \mid wv \in N\}, \\
        &\ell'(v) = \ell(wv).
    \end{align*}
\end{df}
We need to explain how to do corecursion over trees. In order to do so, we are going to use category theory and define recursion at the same time. First we need to define algebras, coalgebras and the morphisms between them.

\begin{df}[Algebra/Coalgebra]
Let \(F\) be an endofunctor of the category \(\textbf{Set}\). An \(F\)-algebra is a pair \((A,\alpha)\) of a set \(A\) and a function \(\alpha : F(A) \longrightarrow A\). An algebra morphism from \((A,\alpha)\) to \((B,\beta)\) is just a function \(f : A \longrightarrow B\) such that the following diagram commutes:
\begin{center}
\begin{tikzcd}
F(A) \arrow[r, "Ff"] \arrow[d, "\alpha"'] & F(B) \arrow[d, "\beta"] \\
A \arrow[r, "f"']                         & B                      
\end{tikzcd}
\end{center}

Similarly, an \(F\)-coalgebra is a pair \((C,\gamma)\) of a set \(C\) and a function \(\gamma : C \longrightarrow F(C)\). A coalgebra morphism from \((C,\gamma)\) to \((D,\delta)\) is a function \(g : C \longrightarrow D\) such that the following diagram commutes:
\begin{center}
\begin{tikzcd}
C \arrow[r, "g"] \arrow[d, "\gamma"'] & D \arrow[d, "\delta"] \\
F(C) \arrow[r, "Fg"']                 & F(D)                 
\end{tikzcd}
\end{center}
\end{df}

An intuition to think about algebras and coalgebras is to imagine the functor \(F\) is representing some structure. Then, the objects of an algebra \(A\) are \emph{created} using the structures of shape \(F\) over \(A\). Conversely, the objects of a coalgebra \(C\) are \emph{destructed} into structures of shape \(F\) over \(C\). Using (co)algebras as objects and (co)algebra morphisms as arrows, we can define a category we will call it \(F\)\textbf{-(Co)Alg}. We say that a (co)algebra is initial (final), in case it is the\footnote{We are justified to talk about the initial (final) (co)algebra, since if it exists it is unique up to (co)algebra isomorphism.} initial (final) object of the corresponding category. Finally, we can define what it means to define a function by recursion or corecursion.

\begin{df}[Recursion/Corecursion]
    Let \((A,\alpha)\) be the initial algebra of an endofunctor \(F\) and \(B\) be a set. We say that \(f : A \longrightarrow B\) has been defined by \emph{recursion} iff there exists a function \(\beta : F(B) \longrightarrow B\) such that \(f\) is the only algebra morphism from \((A,\alpha)\) to \((B,\beta)\).

    Let \((D,\delta)\) be the final coalgebra of an endofunctor \(F\) and \(C\) be a set. We say that \(g : C \longrightarrow D\) has been defined by \emph{corecursion} iff there exists a function \(\gamma : C\longrightarrow F(C)\) such that \(g\) is the only coalgebra morphism from \((C,\gamma)\) to \((D,\delta)\).
\end{df}

We can define an endofunctor of \(\textbf{Set}\), \(\mathcal{T}_L\), such that non-wellfounded (finitely branching) \(L\)-labelled trees are its final coalgebra and finite \(L\)-labelled trees are its initial algebra. 

\begin{df}[Tree endofunctor]
Let \(L\) be a set of labels. We define the \(\textbf{Set}\) endofunctor \(\mathcal{T}_L\) as:
\begin{align*}
    &\mathcal{T}_L(A) = L \times A^*,\\
    &\mathcal{T}_L(f: A \longrightarrow B) = \text{id}_L \times \text{map}_f,
\end{align*}
where \(\text{map}_f : A^* \longrightarrow B^*\) is the pointwise application of \(f\).
\end{df}

We need to find the functions that make the finite trees an initial algebra and the non-wellfounded trees a final coalgebra. These functions are well-known, we call them \(\construct : L \times \tree^*\longrightarrow \tree\) and \(\destruct: \tree \longrightarrow L \times \tree^*\).

\begin{df}
    We define the function \(\construct : L \times \tree^*\longrightarrow \tree\) as \newline \(\construct(a,((N_0,\ell_0),\ldots,(N_{n-1},\ell_{n-1}))) = (N,\ell)\) where:
    \begin{align*}
        &N = \{\epsilon\} \cup \bigcup_{i < n} \{iw \mid w \in N_i \}, \\
        &\ell(w) = \begin{cases}
            a & \text{if }w = \epsilon, \\
            \ell_i(v) & \text{if }w = iv.
        \end{cases}
    \end{align*}
    We define the function \(\destruct: \tree \longrightarrow L \times \tree^*\) as:
    \[
    \destruct(N,\ell) = (\ell(\epsilon),(\textsf{succ}_i(N,\ell))_{i < \arity{T}{\epsilon}})
    \]
    where \(\textsf{succ}_i(T) = \subtree{T}{i}\) for \(i < \arity{T}{\epsilon}\), called the \(i\)-th successor of \(T\).
\end{df}
In simple words, given a label \(a\) and a finite sequence of trees \(T_0,\ldots,T_{n-1}\) we have that \(\construct\) creates the tree whose root has \(a\) as label and \(T_i\) as the \(i\)-th successor. Similarly given a tree \(T\), \(\destruct\) will return a pair with the label of the root and the sequence of trees which are successors of the root, in order. Note that if we restrict the domain to finite trees, we will obtain finite trees in the codomain for both functions. It is easy to check that:

\begin{lem}
    \((\tree,\destruct)\) is the final coalgebra of \(\mathcal{T}_L\) and \((\ftree,\construct)\) is the initial algebra of \(\mathcal{T}_L\).
\end{lem}
In other words, we can define functions to trees by corecursion and from finite trees by recursion.

\section{Non-wellfounded and cyclic proofs}
\label{sec:org7aa93a8}

We work with formulas in the language described by the following BNF:
\[
\phi ::= p \mid \bot \mid \phi \to \phi \mid \phi \wedge \phi \mid \phi \vee \phi \mid \Box \phi,
\]
where \(p\) is a propositional variable. Note that, since our base logic is intuitionistic, \(\Diamond\) is not definable with \(\Box\). In other words, we are working in the \(\Box\)-fragment of modal logic.\footnote{\label{X} In particular \(\textup{i}\mathsf{GL}\) will be the smallest modal logic with intuitionistic propositional logic and the non-logical axioms of \(\textsf{GL}\) formulated with \(\Box\). There is an alternative approach to intuitionistic \(\textsf{GL}\) to consider also the diamond intuitionistically, called \(\textup{I}\mathsf{GL}\). Note that even the \(\Box\)-fragment of these logics is not the same, so they must not be identified. The reader interested in \(\textup{I}\mathsf{GL}\) should consult \cite{das_et_al:LIPIcs.CSL.2024.22}.}

A sequent is just a pair \((\Gamma,\phi)\) where \(\Gamma\) is a finite multiset of formulas and \(\phi\) is a formula. We will use \(\sequent\) to denote the set of all sequents. In other words, we work with 2-sided single conclusion sequents. A \emph{rule instance} is a pair consisting in a finite sequence of sequents and a sequence, called premises and conclusion. A \emph{rule} is just a set of rule instances, let \(\rul\) be the set consisting in all rules (i.e. the set with all sets of rule instances). We are interested in the following rules:

\[\prfbyaxiom{Prop}{\Gamma,p\Rightarrow p} \hspace{3.5cm} \prfbyaxiom{Abs}{\Gamma,\bot \Rightarrow \phi}\]
\[
\prftree[r]{\(\wedge\)L}
{\Gamma,\phi,\psi \Rightarrow \chi}
{\Gamma,\phi \wedge \psi \Rightarrow \chi}
\hspace{3cm}
\prftree[r]{\(\wedge\)R}
{\Gamma \Rightarrow \phi}
{\Gamma \Rightarrow \psi}
{\Gamma \Rightarrow \phi \wedge\psi}
\]
\[\hspace{0.5cm}
\prftree[r]{\(\vee\)L}
{\Gamma,\phi \Rightarrow \chi}
{\Gamma,\psi \Rightarrow \chi}
{\Gamma,\phi\vee\psi \Rightarrow \chi}
\hspace{1.5cm}
\prftree[r]{\(\vee\)R\(_1\)}
{\Gamma \Rightarrow \phi}
{\Gamma \Rightarrow \phi \vee \psi}
\prftree[r]{\(\vee\)R\(_2\)}
{\Gamma \Rightarrow \psi}
{\Gamma \Rightarrow \phi \vee \psi}
\]
\[
\prftree[r]{\(\to\)L}
{\Gamma,\phi \to\psi \Rightarrow \phi}
{\Gamma,\psi \Rightarrow \chi}
{\Gamma,\phi \to \psi \Rightarrow \chi}
\hspace{2cm}
\prftree[r]{\(\to\)R}
{\Gamma,\phi \Rightarrow \psi}
{\Gamma\Rightarrow \phi \to \psi}
\hspace{1cm}
\]
\[
\prftree[r]{\(\Box_{\textsf{K4}}\)}
{\Gamma,\Box\Gamma\Rightarrow \phi}
{\Pi,\Box\Gamma \Rightarrow \Box\phi}
\hspace{2cm}
\prftree[r]{\(\Box_{\textsf{GL}}\)}
{\Gamma,\Box\Gamma,\Box \phi\Rightarrow \phi}
{\Pi,\Box\Gamma \Rightarrow \Box\phi}
\]
where \(\Gamma,\Pi\) are multisets of formulas, \(p\) is a propositional variable and \(\phi,\psi,\chi\) are formulas.
The first 9 rules are called propositional rules, while the last 2 are called modal rules.

From now on, we assume that with tree we mean \(\text{Seq}\times\text{Rul}\)-labelled tree. This permits us to talk about the premises, conclusion and rule of any node \(w\) as follows.
\begin{df}
    Let \(\pi\) be a tree and \(w\) one of its nodes. We define:
    \begin{align*}
        &\pi\text{-prem}(w) = (\textsf{fst}(\ell(wi)))_{i < \arity{\pi}{w}}, \\
        &\pi\text{-concl}(w) = \textsf{fst}(\ell(w)), \\
        &\pi\text{-rule}(w) = \textsf{snd}(\ell(w)),
    \end{align*}
    where \(\textsf{fst}\) is the first projection from an ordered pair and \(\textsf{snd}\) the second projection.
\end{df}

A proof in \(\textup{i}\mathsf{GL}\) is just a standard finite prooftree generated by the propositional logic rules and the modal rule \(\Box_{\textsf{GL}}\). Let us define non-wellfounded and cyclic proofs in \(\textup{i}\mathsf{K4}\).

\begin{df}
A \emph{non-wellfounded proof in \(\textup{i}\mathsf{K4}\)} is a tree \(\pi\) such that:

\begin{enumerate}
\item For any node \(w\), \((\pi\text{-prem}(w),\pi\text{-concl}(w))\) is an instance of the rule \(\pi\text{-rule}(w)\) and \(\pi\text{-rule}(w)\) is either a propositional rule or \(\Box_{\textsf{K4}}\).
\item (Progress condition) For any branch \(w\), there are infinitely many \(i\)'s with \(\pi\text{-rule}(w{\upharpoonright} i) = \Box_{\textsf{K}4}\).
\end{enumerate}
We will write \(\vdash_{\textup{i}\mathsf{K4}_{\infty}} S\) to mean that \(\pi\) is a non-wellfounded proof in \(\textup{i}\mathsf{K4}\) and \(\pi\text{-concl}(\epsilon) = S\). Also, we will denote the collection of non-wellfounded proofs in \(\textup{i}\mathsf{K4}\) as \(\pf{\textup{i}\mathsf{K4}_\infty}\).
\end{df}

\begin{df}
A \emph{cyclic proof in \(\textup{i}\mathsf{K4}\)} is a pair \((\tau,b)\) such that:

\begin{enumerate}
\item \(\tau\) is a finite tree and \(b\) is a partial function from leaves of \(\tau\) to the internal nodes of \(\tau\) called the \emph{backlink function}.
\item For any node \(w \not\in\text{dom}(b)\), we have that \((\pi\text{-prem}(w),\pi\text{-concl}(w))\) is an instance of the rule \(\pi\text{-rule}(w)\) and \(\pi\text{-rule}(w)\) is a propositional rule or \(\Box_{\textsf{K4}}\).
\item For any node \(w \in \text{dom}(b)\), we have that:
\begin{enumerate}
\item \(b(w)\) is a (strict) initial segment of \(w\).
\item \(\pi\text{-concl}(w) = \pi\text{-concl}(b(w))\) and \(\pi\text{-rule}(w) = \varnothing\).
\item (Progress condition) There is a \(v\) between \(b(w)\) and \(w\) (\(b(w)\) initial segment of \(v\) and \(v\) initial segment of \(w\)), with \(\pi\text{-rule}(v) = \Box_{\textsf{K4}}\).
\end{enumerate}
\end{enumerate}
 We will write \(\vdash_{\textup{i}\mathsf{K4}_{\circ}} S\) to mean that \(\pi\) is a cyclic proof in \(\textup{i}\mathsf{K4}\) and \(\pi\text{-concl}(\epsilon) = S\). Also, we will denote the collection of cyclic proofs in \(\textup{i}\mathsf{K4}\) as \(\pf{\textup{i}\mathsf{K4}_\circ}\).
\end{df}

\subsection{On the need of progress}
\label{sec:progress}

Let us show that with the sequent rules we have chosen we explicitly need a progress condition. For that simply consider the following cyclic proof:
\[
\prftree[r]{\(\to\)R}
{\prfbyaxiom{Cycle to (i)}{ p \to q, q \to p \Rightarrow p}}
{
    \prftree[r]{\(\to\)R}
    {\prfbyaxiom{Prop}{q,q \to p \Rightarrow q}}
    {\prfbyaxiom{Prop}{q,p \Rightarrow p}}
    {q, q \to p \Rightarrow p}
}
{\text{(i)} \hspace{.5cm} p \to q, q \to p \Rightarrow p}
\]
Clearly \(p \to q, q \to p \Rightarrow p\) should not be a provable sequent. Since this cyclic proof can also be seen as an infinitary proof it follows that the progress condition is needed in both systems.

\section{Infinitary Proof Translation}
\label{sec:orgfb3ff73}

Thanks to the definition of corecursion we will be able to define a proof translation from a function \(\alpha : \ftree \longrightarrow (\text{Seq} \times \text{Rul}) \times (\ftree)^*\). However, not any function of that shape will give a function from proofs to proofs, in the following definition we enumerate the necessary conditions for this to happen.

\begin{df}[Infinitary Proof Translation]
Let \(\alpha : \ftree \longrightarrow (\text{Seq} \times \text{Rul}) \times (\ftree)^*\). We say that it is an infinitary proof translation iff for any \(\pi \in \pf{\textup{i}\mathsf{GL}}\), if we denote:
\begin{align*}
&\alpha(\pi) = ((S,R),(\tau_0,\ldots, \tau_{n-1})), \\
&\alpha(\tau_i) = ((S_i, R_i),\ldots),
\end{align*}
then the following conditions are satisfied:
\begin{enumerate}
\item \(\tau_0,\ldots, \tau_{n-1} \in \pf{\textup{i}\mathsf{GL}}\).
\item The following is a rule instance of \(\textup{i}\mathsf{K4}\):
    \[
        \prftree[r]{$R$}{S_0}{\cdots}{S_{n-1}}{S}
    \]
\item If \(R \neq \Box_{\textsf{K4}}\), then \(\text{height}(\tau_0),\ldots,\text{height}(\tau_{n-1}) < \text{height}(\pi)\).
\end{enumerate}
Given such \(\alpha\) we define \(\text{trans}_\alpha\) as the only coalgebra morphism from \((\ftree,\alpha)\) to \((\tree,\destruct)\). This implies that \[\text{trans}_{\alpha} = \construct \circ (\text{id} \times \textsf{map}_{\text{trans}_{\alpha}}) \circ \alpha.\]
\end{df}

We want to show that if \(\alpha\) is an infinitary proof translation and \(\pi\) is a (finite) proof in \(\textup{i}\mathsf{GL}\), then \(\text{trans}_{\beta}(\pi)\) is a non-wellfounded proof in \(\textup{i}\mathsf{K4}\). First we need the following technical lemma:

\begin{lem}
\label{lm:story-of-creation-tree}
Let \(\alpha\) be an infinitary proof translation and \(\pi \in \pf{\textup{i}\mathsf{GL}}\). If \(w\) is a node of \(\text{trans}_{\alpha}(\pi)\), then there is an unique sequence of finite \(\textup{i}\mathsf{GL}\)-proofs, \((\iota_{i})_{i \leq \text{length}(w)}\), such that \(\iota_0 = \pi\), \(\subtree{\text{trans}_{\alpha}(\pi)}{w{\upharpoonright i}} = \text{trans}_{\alpha}(\iota_i)\) for \(i \leq\text{length}(w)\) and \(\iota_{i+1} = \textsf{succ}_{w_i}(\alpha(\iota_i))\) for \(i <\text{length}(w)\).

Similarly, if \(w\) is a branch of \(\text{trans}_{\alpha}(\pi)\), then there is an unique sequence of finite \(\textup{i}\mathsf{GL}\)-proofs, \((\iota_{i})_{i \in \omega}\), such that \(\iota_0 = \pi\), \(\subtree{\text{trans}_{\alpha}(\pi)}{w{\upharpoonright i}} = \text{trans}_{\alpha}(\iota_i)\) and \(\iota_{i+1} = \textsf{succ}_{w_i}(\alpha(\iota_i))\) for \(i \in \mathbb{N}\).

\end{lem}
\begin{proof}
The result for nodes is proven by induction in the length of the node. To prove the result for a branch \(b\) it is enough take the union of the sequences given by applying the result for nodes to \(b \restricts i\) for each \(i \in \mathbb{N}\).
\end{proof}

With this lemma we can show that infinitary proof translations transform finite proofs in \(\textup{i}\mathsf{GL}\) into infinitary proofs in \(\textup{i}\mathsf{K4}\):

\begin{thm}
\label{thm:inf-proof-trans}
If \(\alpha\) is an infinitary proof translation, then
\[
\text{trans}_{\alpha} : \pf{\textup{i}\mathsf{GL}} \longrightarrow \pf{\textup{i}\mathsf{K4}_\infty}.
\]
\end{thm}
\begin{proof}
    Let \(\iota \in \pf{\textup{i}\mathsf{GL}}\), we have to check that every node in \(\pi = \text{trans}_\alpha(\iota)\) is an instance of a rule of \(\textup{i}\mathsf{K4}\) and that the progress condition is fulfilled.

    Proof that every node is the instance of a rule. Let \(w\) be a node of \(\pi\) of length \(n\). By Lemma \ref{lm:story-of-creation-tree} we have that there is a sequence of finite trees \((\iota_i)_{i \leq n}\). By the statement of the Lemma \(\iota_0 = \iota\) and \( \iota_{i+1} = \textsf{succ}_{w_i}(\alpha(\iota_i))\). This together with the first condition of proof translation gives that for each \(i\), \(\iota_i \in \pf{\textup{i}\mathsf{GL}}\). Using Lemma \ref{lm:story-of-creation-tree} with \(wj\) for \(j < \arity{\pi}{w}\) we get sequences \((\iota^j_i)_{i \leq n +1}\) such that \(\iota^j_i = \iota_i\) for \(i \leq n\) and \(\iota^j_{n+1} = \textsf{succ}_j(\alpha(\iota_n))\). We note that since \(\subtree{\text{trans}_\alpha(\iota)}{w} = \text{trans}_\alpha{\iota_n}\) and \(\subtree{\text{trans}_\alpha(\iota)}{wj} = \text{trans}_\alpha{\iota^j_{n+1}}\) the conclusion of rule instance at \(w\) in \(\pi\) is just the result of looking at the sequent given by applying \(\alpha\) tp \(\iota_n\) and the premises to is the sequent given by applying \(\alpha\) to each \(\iota^j_n\).
    Using this and the second condition of infinitary proof translation to \(\iota_n\), since \(\iota_n \in \pf{\textup{i}\mathsf{GL}}\). we get that the node is the instance of a rule, as desired.

    Proof of progress condition. Imagine the result does not fulfill the progress condition, then there is an infinite branch \(b\) in \(\pi\) such that from some point do not contain applications of \(\Box_{\textsf{K4}}\). By applying Lemma \ref{lm:story-of-creation-tree} to this branch \(b\) we would obtain an infinite sequences of finite trees \((\iota_n)_{n \in \mathbb{N}}\). By the third condition of infinitary proof translation we would obtain that from some point the height of these trees is strictly decreasing, absurd.
\end{proof}

\section{From finitary to cyclic, and vice versa}
\label{sec:org9e64476}

In order to define the translation from the finite acyclic system into the cyclic system first we need to obtain non-wellfounded proofs of certain shape. For this, we need the admissibility of the following two rules in \(\textup{i}\mathsf{GL}\):
\[
\prftree[r]{Contract}{\Gamma,\psi,\psi \Rightarrow \phi}{\Gamma,\psi \Rightarrow \phi} \hspace{0.5cm}
\prftree[r]{L\"ob}{\Gamma,\Box\Gamma,\Box\phi \Rightarrow \phi}{\Gamma,\Box\Gamma \Rightarrow \phi}
\]
The admissibility of these rules is proven in \cite{Iris}. If \(\Gamma\) is a finite multiset it can be split into a set (a multiset where all its elements appear exactly once) and a multiset in a unique way. \(\Gamma^s\) will be the resulting set of this splitting and \(\Gamma^m\) will be the multiset.

Thanks to admissibility of contraction we can define a function that, given an \(\textup{i}\mathsf{GL}\)-proof  \(\pi\) of \(\Gamma,\Box \Gamma \Rightarrow \phi\), returns a proof \(\textsf{contract}(\pi)\) of \(\Gamma^s, \Box \Gamma^s \Rightarrow \phi\). Similarly, thanks to admissibility of \(\text{L\"ob}\) we can define a function that given an \(\textup{i}\mathsf{GL}\)-proof \(\pi\) of \(\Gamma,\Box\Gamma,\Box\phi\Rightarrow\phi\), returns a proof \(\textsf{l\"ob}(\pi)\) of \(\Gamma,\Box\Gamma\Rightarrow\phi\).

Using these two functions and the theorem of proof translations, we have the following result:

\begin{thm}
\label{th:wanted-translation}
There is a unique \(h : \pf{\textup{i}\mathsf{GL}} \longrightarrow \pf{\textup{i}\mathsf{K4}_\infty}\), such that \(h\) is the identity in initial sequents, commutes with the logical rules and:
\[
\prftree[noline]{\pi}{\prftree[r]{$\Box_{\textsf{GL}}$}{\Gamma,\Box\Gamma,\Box\phi \Rightarrow \phi}{\Pi,\Box\Gamma\Rightarrow \Box\phi}}\longmapsto 
\prftree[noline]{h(\textsf{contract}(\textsf{l\"ob}(\pi))}{
    \prftree[r]{\(\Box_{\textsf{K4}}\)}
    {\Gamma^s,\Box\Gamma^s \Rightarrow \phi}
    {\Pi,\Box\Gamma^m,\Box\Gamma^s\Rightarrow \Box\phi}
}
\]
In addition, if \(\pi\) is a proof of \(\Gamma \Rightarrow \phi\), then \(h(\pi)\) is also a proof of \(\Gamma \Rightarrow \phi\).
\end{thm}
\begin{proof}
Simple application of \ref{thm:inf-proof-trans}.
\end{proof}

Finally, thanks to this proof translation we can show that any finite acyclic proof can be transformed into a cyclic proof.

\begin{cor}
\label{cr:main-result}
Let \(S\) be a sequent. If \(\vdash_{\textup{i}\mathsf{GL}} S\), then \(\vdash_{\textup{i}\mathsf{K4}_\circ} S\).
\end{cor}
\begin{proof}
This result is based in two observations. First, that thanks to the shape of the rules we have the subformula property. Second, that any branch must have infinitely many applications of the \(\Box_{\textsf{K4}}\) rule, but choosing the translation we defined before the premise of such an application is determined by a set. Using these two facts together we can get the required repetition of sequents to define a cyclic proof.
\end{proof}

The argument to obtain a proof from the cyclic system to the finite acyclic system is totally analogous to the classical case proven in \cite{ShamkanovGl}. We can conclude the desired result:
\begin{thm}
Let \(S\) be a sequent. Then \(\vdash_{\textup{i}\mathsf{GL}}S\) iff \(\vdash_{\textup{i}\mathsf{K4}^\circ}S\).
\end{thm}

\section{Conclusion}

We have provided a proof of the equivalence between a finite acyclic system with rules of \(\textup{i}\mathsf{GL}\) and a (finite) cyclic system with rules of \(\textup{i}\mathsf{K4}\) using a corecursive translation of proofs. In order to do this, we exploited that while performing a corecursion we can use the admissible rules in the finite acyclic system to obtain non-wellfounded proofs of the desired shape. 

\bibliographystyle{eptcs}
\bibliography{bib}
\end{document}